\documentclass[a4paper, 11pt]{amsart}
\usepackage[latin1]{inputenc}
\usepackage[ T1]{fontenc}
\usepackage[english]{babel}
\usepackage{amssymb}
\usepackage{amsmath}
\usepackage{amsthm}
\usepackage{amscd}
\usepackage{amsfonts}
\usepackage{stmaryrd}
\usepackage{pb-diagram}
\usepackage{epic,eepic,epsfig}
\usepackage{a4wide}
\usepackage{nextpage}
\usepackage{fancyhdr}
\usepackage{enumerate}
\usepackage{hyperref}
\usepackage{float}
\usepackage{comment}

\newtheorem{prop}{Proposition}[section]
\newtheorem{corollary}[prop]{Corollary}

\newtheorem{lemma}[prop]{Lemma}
\newtheorem{thm}[prop]{Theorem}

\newtheorem{theorem}[prop]{Theorem}

\theoremstyle{definition}
\newtheorem{rem}[prop]{Remark}

\newtheorem{definition}[prop]{Definition}

\newcommand{\Ima} {\mathop{\mathrm{Im}}}

\newcommand{\degr} {\mathop{\mathrm{deg}}}
\newcommand{\degrar} {\mathop{\widehat{\mathrm{deg}}}}

\newcommand{\Spec} {\mathop{\mathrm{Spec}}}

\newcommand{\SL} {\mathrm{SL}}

\newcommand{\eps}{\varepsilon}
\newcommand{\psitil}{\widetilde{\psi}}

\renewcommand{\Im}{\mathop{\mathrm{Im}}}
\renewcommand{\Re}{\mathop{\mathrm{Re}}}

\def\BC{\mathbb{C}}

\def\BZ{\mathbb{Z}}
\def\BH{\mathbb{H}}

\def\CF{\mathcal{F}}

\usepackage[OT2,T1]{fontenc}
\DeclareSymbolFont{cyrletters}{OT2}{wncyr}{m}{n}
\DeclareMathSymbol{\Sha}{\mathalpha}{cyrletters}{"58}

\begin{document}

\title{Explicit bounds on the coefficients of modular polynomials and the size of $X_0(N)$}
\author{}
\author{Florian Breuer, Desir\'ee Gij\'on G\'omez, Fabien Pazuki}

\address{School of Information and Physical Sciences, The University of Newcastle,
University Drive, Callaghan, NSW 2308, Australia.}
\email{Florian.Breuer@newcastle.edu.au}

\address{Department of Mathematical Sciences, University of Copenhagen,
Universitetsparken 5, 
2100 Copenhagen \O, Denmark.}
\email{dgg@math.ku.dk}

\address{Department of Mathematical Sciences, University of Copenhagen,
Universitetsparken 5, 
2100 Copenhagen \O, Denmark, and Universit\'e de Bordeaux, 33405 Talence, France.}
\email{fpazuki@math.ku.dk}

\thanks{The authors thank Pascal Autissier, Joe Silverman, and Emmanuel Ullmo, for conversations around this topic at the occasion of the Hindry 65 conference in Bordeaux. They also thank Autissier for comments on an earlier version of the text. They thank Riccardo Pengo and Paolo Dolce for providing the reference \cite{DM}. The authors were supported by the IRN GandA (CNRS). The third author is supported by ANR-20-CE40-0003 Jinvariant.}

\thispagestyle{empty}

\maketitle

\noindent \textbf{Abstract.}
We give explicit upper and lower bounds on the size of the coefficients of the modular polynomials $\Phi_N$ for the elliptic $j$-function. These bounds make explicit the best previously known asymptotic bounds. We then give an explicit version of Silverman's Hecke points estimates. Finally, we give an asymptotic comparison between the Faltings height of the modular curve $X_0(N)$ and the height of the modular polynomial $\Phi_N$.

{\flushleft
\textbf{Keywords:} Modular polynomials, modular curves, elliptic curves, heights.\\
\textbf{Mathematics Subject Classification:} 11F32, 11G05, 11G50, 14G40. }
\begin{center}
---------
\end{center}

\section{Introduction}

Modular curves play a central role in modern arithmetic questions. They are a key feature in the solution of famous diophantine equations, in the study of the Mordell-Weil group of elliptic curves (both for the torsion subgroup and for the Birch and Swinnerton-Dyer conjecture) and in isogeny-based cryptography. It is thus useful to be able to represent explicitly these curves and to estimate how complicated their models are.

A classical way to estimate complexity of models is via height theory. For any non-zero polynomial $P$ in one or more variables and complex coefficients we define its {\em height} to be 
\[
h(P) := \log \max|c|, \quad\text{where $c$ ranges over all coefficients of $P$.}
\]

Let $N$ be a positive integer and denote by $\Phi_N = \Phi_N(X,Y) \in \BZ[X,Y]$
the modular polynomial for the elliptic $j$-function. It vanishes at pairs of $j$-invariants of elliptic curves linked by a cyclic $N$-isogeny, see \cite[Chapter 5]{Lang}. The equation $\Phi_N(X,Y)=0$ is a plane affine integral model for the modular curve $X_0(N)$ (but not in general a smooth model).

Paula Cohen Tretkoff \cite{Coh} proved that when $N$ tends to $+\infty$
\begin{equation}\label{Paula}
   h(\Phi_N) = 6\psi(N)\big[\log N - 2\kappa_N + O(1)\big], 
\end{equation}
where
\[
\psi(N) = N\prod_{p|N}\left(1 + \frac{1}{p}\right)\quad\text{and}\quad
\kappa_N = \sum_{p|N}\frac{\log p}{p}.
\]

Work of Autissier \cite{Aut} and Breuer-Pazuki \cite{BP} show that one may profitably replace $\kappa_N$ with $\lambda_N$, where \[\lambda_N := \sum_{p^n\|N}\frac{p^n-1}{p^{n-1}(p^2-1)}\log p.\]

The terms $\kappa_N$ and $\lambda_N$ are compared in \cite{BP}, which leads in particular to

\begin{equation}\label{Paula+BP}
   h(\Phi_N) = 6\psi(N)\big[\log N - 2\lambda_N + O(1)\big].
\end{equation}
Numerical computations as reported in \cite{BP} suggest that the bounded term implied by the $O(1)$ in (\ref{Paula+BP}) is smaller than the one in (\ref{Paula}).

As modular polynomials have various cryptographic or algorithmic applications, it is useful to obtain explicit bounds on the $O(1)$ term. In \cite{BrSu}, Br\"oker and Sutherland obtained asymptotically optimal bounds in the case where $N$ is prime, and Pazuki \cite{Paz} provided explicit bounds in the case of general $N$, but these were not quite asymptotically optimal.

The most recent work providing an explicit upper bound is \cite{BP}, where the first and third authors proved that for any $N\geq2$,
\begin{equation}\label{eq:BPold}
    h(\Phi_N) \leq 6\psi(N)\big[\log N -2\lambda_N + \log\log N + 4.436\big].
\end{equation}

The term $\log\log N$ was superfluous, an unfortunate artifact of the method used in \cite{BP}. A natural idea to try to remove it is via equidistribution results. However, that would be at the cost of losing the explicit nature of the upper bound, hence jeopardizing our other efforts. We are nevertheless now able to remove the $\log\log N$ in the following theorem, where we provide both explicit upper and lower bounds.

\begin{thm}\label{thm:main1}
    Let $N \geq 1$. The height of the modular polynomial $\Phi_N(X,Y)$ is bounded by
    \[
    6\psi(N)\big[\log N -2\lambda_N -0.0351\big] \leq h(\Phi_N) \leq 6\psi(N)\big[\log N -2\lambda_N + 9.5387\big].
    \]
\end{thm}

The main new idea to improve the upper bound comes from technical inequalities involving Farey sequences. We use both reduced and non-reduced elements in the upper half plane in the key equation (\ref{SNtau}), which help us obtain better estimates of the Mahler measures at play.

To obtain the lower bound, we use a specialization trick to reduce the calculations to the Mahler measure of the one-variable polynomial $\Phi_N(X,0)=\Phi(X,j(\rho))$, where we can use explicit complex multiplication properties. 

After presenting some preliminaries in Section \ref{sec:prelim}, we prove the upper bound of Theorem \ref{thm:main1} in Section \ref{proof of upper}. We prove the lower bound of Theorem \ref{thm:main1} in Section \ref{proof of lower}.
\\

 As a corollary to Theorem \ref{thm:main1}, we add the following explicit result on Hecke points in Section \ref{Hecke}, giving an explicit version of a result of Silverman \cite{Sil}. For any elliptic curve $E$ defined over $\overline{\mathbb{Q}}$, for any $N\geq2$ and any cyclic subgroup $C\subset E(\overline{\mathbb{Q}})$ of order $N$, denote by $j_{E/C}$ the $j$-invariant of the isogenous elliptic curve $E/C$.

\begin{theorem}\label{thm:explicit Silverman}
    Let $E$ be an elliptic curve defined over $\overline{\mathbb{Q}}$ with $j$-invariant $j_E$. Let $h_\infty(\cdot)$ denote the absolute logarithmic Weil height. For any $N\geq 2$, one has 
    \begin{enumerate}
    \item[(a)] $\displaystyle h_\infty(j_E)-\frac{1}{\psi(N)}\sum_{\substack{C\; cyclic \\ \# C=N}} h_\infty(j_{E/C})$ \\
    $\displaystyle\geq -\frac{h(\Phi_N)}{\psi(N)}  - \frac{2\log(\psi(N) + 1)}{\psi(N)} \geq -6\log N +12\lambda_N - 58.34.$
    \item[(b)] $\displaystyle h_\infty(j_E)-\frac{1}{\psi(N)}\sum_{\substack{C\; cyclic \\ \# C=N}} h_\infty(j_{E/C})$ \\
    $\displaystyle    \leq 6.67 + 6\min\big\{0, \log\big(1 + h_\infty(j_E)\big) - \log N + 2\lambda_N + 0.25\big\}.$
    \end{enumerate}
\end{theorem}

The proof is given in Section \ref{Hecke}. It combines Silverman's method, Mahler measure estimates, and the explicit bounds from Theorem \ref{thm:main1}. 

The height of $\Phi_N$ is a way to measure the size of the curve $X_0(N)$. But there are other ways of measuring the size of $X_0(N)$: the Faltings height of the curve, the Faltings height of its jacobian $J_0(N)$, the height of a Hecke correspondence with respect to a carefully chosen metrized line bundle, the auto-intersection of the Arakelov canonical sheaf, are all used in the literature. One could even think of the size of classical Heegner points on the modular jacobian as a way to measure the complexity of $J_0(N)$, hence of $X_0(N)$. So what is the size of $X_0(N)$? We gather in the following theorem some asymptotic results that are easy to derive from the existing literature, and which explain that the height of $\Phi_N$, despite being elementary, captures some of this deeper information. 

\begin{theorem}\label{what is the size}
    We have the following properties.
    \begin{enumerate}
        \item[(a)] Let $h_{\mathrm{Falt}}$ denote the stable Faltings height as recalled in Definition \ref{faltings}. For any integer $N\geq1$, one has the equality $h_{\mathrm{Falt}}(X_0(N))=h_{\mathrm{Falt}}(J_0(N))$. Then when $N$ is square-free and coprime to $6$ and tends to infinity, one has $$h_{\mathrm{Falt}}(X_0(N))\sim \frac{1}{6^3} h(\Phi_N).$$
        \item[(b)] Let $T_N$ be the Hecke correspondence in $\mathbb{P}^1\times \mathbb{P}^1$ and let $\hat{\mathcal{L}}$ be the associated metrized line bundle as given by Autissier in \cite{Aut}. Then when $N$ tends to infinity one has $$ h_{\hat{\mathcal{L}}}(T_N) \sim 2 h(\Phi_N).$$
        \item[(c)] Let $k$ be a quadratic field of discriminant $D_k$, of class number $h_k$, with $2u_k$ roots of unity. Assume $D_k<0$, $D_k=1(\mathrm{mod}\, 4)$, and consider $x_{D_k}\in{X_0(N)}$ the related Heegner point for each compatible $N$. It gives rise to a cycle $c_{D_k}=(x_{D_k})-(\infty)\in{J_0(N)}$. Then when $N$ tends to infinity, one has $$\hat{h}_{J_0(N)}(c_{D_k})\sim \frac{h_k u_k}{6\psi(N)}h(\Phi_N).$$
        \item[(d)] Let $\overline{\omega}^2$ denote the auto-intersection of the Arakelov canonical sheaf of the minimal regular model of $X_0(N)$, for any $N\geq2$ coprime to $6$. Then when $N$ tends to infinity, one has $$\overline{\omega}^2\sim \frac{1}{24}h(\Phi_N).$$
    \end{enumerate}
\end{theorem}

We prove Theorem \ref{what is the size} in Section \ref{what?}.

\section{Preliminaries}\label{sec:prelim}

Denote by $\BH = \{ \tau \in\BC \;:\; \mathrm{Im}(\tau) > 0 \}$ the upper half-plane, on which $\SL_2(\BZ)$ acts via fractional linear transformations. A fundamental domain for this action is 
\[
\CF=\left\{\tau\in\BH\; : \; \vert \tau\vert\geq 1, \;  -\frac{1}{2} < \Re\tau \leq \frac{1}{2} \; \mathrm{and}\;
\Re\tau \geq 0 \; \mathrm{if}\; |\tau|=1\right\}.
\]

For any $\tau\in\BH$, we denote by $\tilde{\tau}\in\CF$ the unique representative in this fundamental domain of the $\SL_2(\BZ)$-orbit of $\tau$.

The $j$-function $j : \BH \to \BC$ is $\SL_2(\BZ)$-invariant, and satisfies a $q$-expansion of the form
\[
j(\tau) = q^{-1} + 744 + 196884q+\cdots, \quad\text{where}\quad q = e^{2\pi i\tau}.
\]
We will also consider the modular discriminant function $\Delta : \BH \to \BC$, which is a cusp form of weight 12 for $\SL_2(\BZ)$, and we choose to normalise it such that its $q$-expansion is
\begin{equation}\label{deltadef}
\Delta(\tau) = q\prod_{n=1}^{\infty}(1-q^n)^{24} = q -24q^2 + 252q^3 + \cdots.
\end{equation}

This modular form plays a key role in this paper. Let us start by computing in the next lemma two special values which will be used in the sequel.

\begin{lemma}\label{deltarho}
Let $\Delta$ be the discriminant modular form, normalized as in (\ref{deltadef}). We have
\begin{enumerate}
    \item[(a)] 
    $\displaystyle\Delta(\rho)=-\frac{3^3}{(2\pi)^{24}}\Gamma\left(\frac{1}{3}\right)^{36},\quad$ where $\rho = e^{\frac{i\pi}{3}}$ and $\Gamma$ stands for Euler's Gamma function, and
    \item[(b)] 
    $\displaystyle \Delta(i)=\frac{1}{2^{24}\pi^{18}}\Gamma\left(\frac{1}{4}\right)^{24}.$
    \noindent 
\end{enumerate}
\end{lemma}
\begin{proof}
    Let us start with (a). We have classically $(2\pi)^{12}\Delta(\rho)=g_2(\rho)^3-27g_3(\rho)^2$, with $g_2, g_3$ the normalized Eisenstein series, and $g_2(\rho)=0$ is a direct computation. For the value $g_3(\rho)$, we work with the elliptic curve in complex Weierstrass form $y^2=4x^3-4$, which has period lattice $\Lambda=\omega\mathbb{Z}+\rho \omega\mathbb{Z}$, with period $$\omega=2\int_{1}^{+\infty}\frac{dt}{\sqrt{4t^3-4}}=\int_{1}^{+\infty}\frac{dt}{\sqrt{t^3-1}}=\frac{1}{3}B\left(\frac{1}{6},\frac{1}{2}\right)=\frac{\Gamma(\frac{1}{3})^3}{2^{\frac{4}{3}}\pi},$$ where $B(.,.)$ is the Euler $B$ function, as classically defined for any complex numbers $z_1, z_2$ with positive real part by $$B(z_1,z_2):=\int_0^1t^{z_1-1}(1-t)^{z_2-1}dt=\frac{\Gamma(z_1)\Gamma(z_2)}{\Gamma(z_1+z_2)}.$$ Writing in generic Weierstrass form $4x^3-4=4x^3-g_2x-g_3$, we simply read off $g_2(\Lambda)=0$ and $g_3(\Lambda)=4$. We can now compute $g_3(\Lambda)=\omega^{-6}g_3(\mathbb{Z}+\rho\mathbb{Z})=\omega^{-6}g_3(\rho)$, hence $g_3(\rho)^2=4^2\omega^{12}$, which gives the claim for $\Delta(\rho)$.
    
    We treat part (b) similarly: $(2\pi)^{12}\Delta(i)=g_2(i)^3-27g_3(i)^2$, and $g_3(i)=0$ is a direct computation. For the value $g_2(i)$, we work with the elliptic curve in complex Weierstrass form $y^2=4x^3-4x$, which has period lattice $\Lambda=\omega_0\mathbb{Z}+i \omega_0\mathbb{Z}$, with period $$\omega_0=2\int_{1}^{+\infty}\frac{dt}{\sqrt{4t^3-4t}}=\int_{1}^{+\infty}\frac{dt}{\sqrt{t^3-t}}=\frac{1}{2}B\left(\frac{1}{4},\frac{1}{2}\right)=\frac{\Gamma(\frac{1}{4})^2}{2^{\frac{3}{2}}\pi^{\frac{1}{2}}}.$$ Writing in generic Weierstrass form $4x^3-4x=4x^3-g_2x-g_3$, we read off $g_3(\Lambda)=0$ and $g_2(\Lambda)=4$. We can now compute $g_2(\Lambda)=\omega_0^{-4}g_2(\mathbb{Z}+i\mathbb{Z})=\omega_0^{-4}g_2(i)$, hence $g_2(i)^3=4^3\omega_0^{12}$, which gives the claim for $\Delta(i)$.
\end{proof}

\begin{rem}
    From old work of Hurwitz, one can also derive another expression of $\Delta(i)$ using another period. From equation (7) page 201 of \cite{Hur} we get $$\Delta(i)=\frac{2^{18}}{(2\pi)^{12}}\left(\int_0^1\frac{dt}{\sqrt{1-t^4}}\right)^{12}.$$
\end{rem}

Our first analytical tool is the following result, which is a refinement of (3.18) of \cite{Paz}.

\begin{lemma}\label{logmaxbound}
    Let $f(\tau) = \log\max\big\{|\Delta(\tau)|, |j(\tau)\Delta(\tau)|\big\}$. Then for all $\tau\in\mathcal{F}$,
    \[
    -5.5335 <  f(\tau) \leq f(i) =
    \log\left(\frac{3^3}{2^{18}\pi^{18}}\Gamma\left(\frac{1}{4}\right)^{24}\right) < 1.1266.
    \]
\end{lemma}

\begin{proof}
    We have 
    \[
    j(\tau) = \frac{g_2(\tau)^3}{(2\pi)^{12}\Delta(\tau)},
    \]
    where
    $g_2(\tau)$ is again the normalized Eisenstein series of weight 4. 
    Thus
    \[
    f(\tau) = \left\{ 
    \begin{array}{ll} \log|\Delta(\tau)| & \text{if $|j(\tau)\vert < 1$} \\
    3\log|g_2(\tau)| - 12\log(2\pi) & \text{if $|j(\tau)| \geq 1$.}
    \end{array}\right.
    \]
    
    The boundary of $\CF$ consists of a circular arc $C$ from $\rho$ to $\rho^2$, where $\rho=e^{\frac{i\pi}{3}}$, 
    as well as the two
    vertical half-lines $L$ from $\rho$ to $i\infty$ and $L'$ from $\rho^2$ to $i\infty$. 
    
    Since $j(\tau)$ has simple zeros at $\rho$ and $\rho^2$, and no other zeroes near $\CF$, 
    one finds that $|j(\tau)|\leq 1$ in small neighbourhoods of these two points.
    Their intersection with $\CF$ consists of two connected components,
    $D\cup D' = \{\tau \in\CF \;:\; |j(\tau)| \leq 1\}$, where $\rho\in D$ and $\rho^2\in D'$. 

    By the Maximum Modulus Principle, $f(\tau)$ attains its extrema either at the cusp $i\infty$ or on the 
    boundary components $L'$, $C$, $L$, $\partial D$ and $\partial D'$.

    Using SageMath \cite{Sage}, we computed $f$ restricted to these boundary components. The results are symmetric around the imaginary
    axis, so Figure \ref{plotf} shows the plot of $f(\tau)$ for $\tau$ on the contour from $i$ via $\rho$ to $i\infty$, as well as on $\partial D$.

    \begin{figure}
        \centering
        \includegraphics[width=0.3\textwidth]{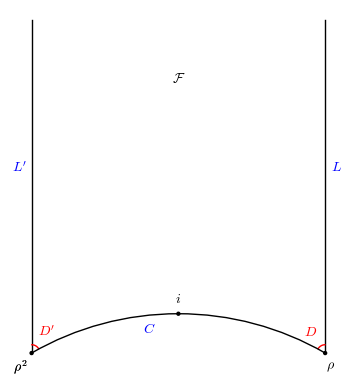}
        \includegraphics[width=0.65\textwidth]{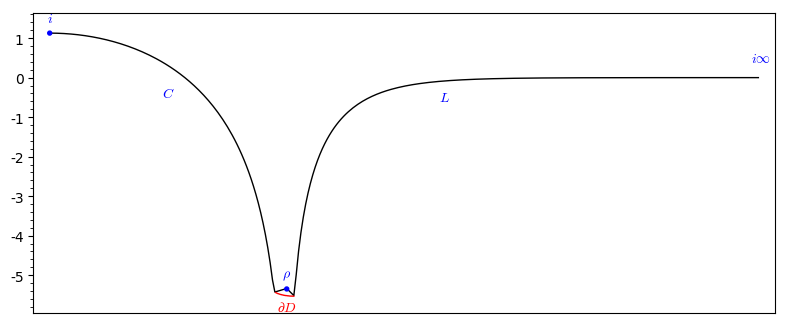}
        \caption{{\em (Left)} The fundamental domain $\CF$. {\em (Right)} Plot of $f(\tau)$ for $\tau\in C\cup L\cup \partial D$ and $\mathrm{Re}(\tau) \geq 0$.}
        \label{plotf}
    \end{figure}

    We find that $f$ attains its maximum at $f(i) < 1.1266$ and its minimum $ > -5.5335$ where $\partial D$ meets $L$. At the cusp, $f(i\infty)=0$, which lies between these two extreme values. The formula for $f(i)$ comes directly from Lemma \ref{deltarho}.

\end{proof}

\begin{rem}
    Repeating the computations in \cite{BP} using the upper bound in Lemma \ref{logmaxbound} instead of \cite[(13)]{BP}, we obtain in the following corollary a slight improvement on the constant in (\ref{eq:BPold}).
    \end{rem}

\begin{corollary}
    Let $N\geq 2$. the height of the modular polynomial $\Phi_N(X,Y)$ is bounded by
    \[
    h(\Phi_N) \leq 6\psi(N)\big[\log N -2\lambda_N + \log\log N + 4.238\big].
    \]
    \qed
\end{corollary}

\section{Proof of the upper bound in Theorem \ref{thm:main1}}\label{proof of upper}

\subsection{Strategy of proof.}

Let us start by denoting, for $N\geq1$, 
\[
C_N=\left\{\left(\begin{matrix}
a & b \\
0 & d
\end{matrix}\right)\; : \; a,b,d\in\BZ, \;
ad=N, \; a\geq1, \; 0\leq b\leq d-1,\; \gcd(a,b,d)=1 \right\}.
\]
We have 
\[
\#C_N = \sum_{d|N}\; \sum_{\substack{0\leq b < d \\ (b,r)=1}} 1
= \sum_{d|N} \frac{d\varphi(r)}{r} 
= \psi(N),
\]
where we denote the gcd $r=(d, \frac{N}{d})$ for each $d$.

The relevance of the matrices in $C_N$ is the following. For each $\gamma\in C_N$ and $\tau\in\BH$,
define 
\[
\tau_\gamma := \gamma(\tau) = \frac{a_\gamma\tau + b_\gamma}{d_\gamma}.
\]
Then the elliptic curves $\BC/\tau\BZ + \BZ$ and $\BC/\tau_\gamma\BZ + \BZ$ are linked by a cyclic isogeny of degree $N$. Conversely, up to isomorphism, all cyclic $N$-isogenies from $\BC/\tau\BZ + \BZ$ are obtained this way. In particular, the modular polynomial $\Phi_N(X,Y)$ satisfies
\[
\Phi_N(X, j(\tau)) = \prod_{\gamma\in C_N}\left(X - j(\tau_\gamma)\right).
\]

An interpolation argument (see Lemma \ref{interpolation}) allows us to estimate the height of $\Phi_N(X,Y)$ in terms of the heights of the specialised polynomials $\Phi_N(X, j(\tau))$ for suitable values of $\tau\in\BH$. These, in turn, are related to their logarithmic Mahler measures:

\[
S_N(\tau) :=
 m\big(\Phi_N(X,j(\tau))\big) =  \sum_{\gamma\in C_N} \log\max\big\{1, |j(\tau_\gamma)|\big\}.
\]

This Mahler measure is now our top priority. Let us work on the formula defining $S_N(\tau)$ and start with equation (14) from \cite{BP}, valid for any $N\geq1$ and $\tau\in{\mathbb{H}}$:
\begin{equation}\label{SNtau}
  S_N(\tau)= \displaystyle{\sum_{\gamma\in C_N}} \log\max\{|\Delta(\tilde{\tau}_\gamma)|,|j(\tau_\gamma)\Delta(\tilde{\tau}_\gamma)|\} 
    +6 \sum_{\gamma\in C_N} \big[\log\Im\tilde{\tau}_\gamma 
    - \log\Im\tau_\gamma\big]
    - \psi(N)\log|\Delta(\tau)|.
      \end{equation}
Recall that here $\tilde{\tau}_\gamma \in\CF$ denotes the representative of $\tau_\gamma$ in the fundamental domain $\CF$. We invoke \cite[Lemme 2.3]{Aut}:
\begin{equation}\label{eq:Aut sum log d/a}
  \sum_{\gamma\in C_N}\log\frac{d_\gamma}{a_\gamma} = \psi(N)(\log N - 2\lambda_N),
\end{equation}
which combined with 
\[
\Im\tau_\gamma = \Im\left( \frac{a_\gamma\tau + b_\gamma}{d_\gamma} \right)
= \frac{a_\gamma}{d_\gamma}\Im\tau
\]
gives
\begin{equation} 
  \label{eq:Aut2.3}
  -\sum_{\gamma\in C_N}\log\Im\tau_\gamma 
  = \psi(N)\big(\log N - 2\lambda_N   -\log\Im\tau\big).
\end{equation}

Inject equality (\ref{eq:Aut2.3}) in equation (\ref{SNtau}) and use the upper bound from Lemma \ref{logmaxbound} (note that $j(\tau_\gamma)=j(\tilde{\tau}_\gamma)$) to get:
\begin{align}
    S_N(\tau) = & \sum_{\gamma\in C_N} \log\max\{|\Delta(\tilde{\tau}_\gamma)|,|j(\tilde\tau_\gamma)\Delta(\tilde{\tau}_\gamma)|\}
     + 6\psi(N)\big[\log N - 2\lambda_N\big] \label{eq:SN0}   \\ \nonumber
    & + 6\sum_{\gamma\in C_N} \log\Im\tilde{\tau}_\gamma - \psi(N)\log\big[|\Delta(\tau)|(\Im\tau)^6\big] \\ 
   \leq & \; 6\psi(N)\big[\log N - 2\lambda_N +  0.1878\big] 
     + 6\sum_{\gamma\in C_N} \log\Im\tilde{\tau}_\gamma 
    - \psi(N)\log\big[|\Delta(\tau)|(\Im\tau)^6\big]. \label{eq:SN}
\end{align}


Our strategy is to set 
$\tau = iy$ with $y\geq 1$
and obtain an explicit upper bound for the sum $\displaystyle{\sum_{\gamma\in C_N} \log\Im\tilde{\tau}_\gamma}$, which we will decompose into a sum with large $d$ and a sum with small $d$:

\begin{equation}
\label{eq:target}
\sum_{\gamma\in C_N} \log\Im\tilde{\tau}_\gamma
= \sum_{\substack{\gamma\in C_N \\ d_\gamma \geq \sqrt{Ny}}} \log\Im\tilde{\tau}_\gamma
+ \sum_{\substack{\gamma\in C_N \\ d_\gamma < \sqrt{Ny}}} \log\Im\tilde{\tau}_\gamma.
\end{equation}

Our strategy is inspired by \cite{Coh}, where the author uses a similar decomposition.

\subsection{Large $d$}

Consider $\gamma\in C_N$ for which $d = d_\gamma \geq \sqrt{Ny}$.
As in \cite{Coh}, we will approximate $\tilde{\tau}_\gamma$ with a representative $\hat{\tau}_\gamma\in \SL_2(\BZ)\tau_\gamma$ satisfying $\Im\hat{\tau}_\gamma \geq \frac{1}{2}$.

We start with the following lemma, which relies on the Farey sequence of order $M$.

\begin{lemma}\label{Farey}
    Let $M\geq1$ be an integer. Then one can express the interval $$I_M=\left[\frac{1}{M+1},\frac{M+2}{M+1}\right)=\bigcup_{k=1}^{M}\bigcup_{\substack{h=1 \\ (h,k)=1}}^{k} I_M\left(\frac{h}{k}\right),$$ as a disjoint union of intervals $I_M \left(\frac{h}{k}\right)$ of the form $[\rho_1, \rho_2)$ containing $\frac{h}{k}$ and such that \[\frac{1}{2Mk}\leq \frac{h}{k}-\rho_1\leq \frac{1}{(M+1)k},\] 
    \[\frac{1}{2Mk}\leq \rho_2-\frac{h}{k}\leq \frac{1}{(M+1)k}.
    \]
\end{lemma}
\begin{proof}
This is \cite[Lemma 3]{Coh}.
\end{proof}
Recall that $\tau = iy$ with $y\geq 1$ and $d \geq \sqrt{Ny}$. Then we set
\[
M := \left\lfloor \frac{d}{\sqrt{Ny}} \right\rfloor \geq 1.
\]

Let 
$\displaystyle\gamma = \left(
\begin{matrix}
    a & b \\
    0 & d
\end{matrix}
\right)\in C_N$. 
If $\frac{b}{d} \in [0, \frac{1}{M+1})$ then we replace $b$ by $b+d$; this merely has the effect of replacing $\tau_\gamma$ by $\tau_\gamma + 1$, which is in the same $\SL_2(\BZ)$-orbit.

Next, choose a matrix
$\delta = \left(\begin{matrix}
s & r \\
k & -h
\end{matrix}\right) \in \SL_2(\BZ)$ 
for which $\frac{b}{d} \in I_M(\frac{h}{k})$ and define 
\[
\hat\tau_\gamma := \delta(\gamma(\tau)).
\]
The entries $s$ and $r$ may be chosen in such a way (multiplying $\delta$ by a suitable translation matrix) that $-\frac{1}{2} < \mathrm{Re}(\hat\tau_\gamma) \leq \frac{1}{2}$.

\begin{lemma}\label{tauhat}
The elements $\hat{\tau}_\gamma$ constructed above satisfy the following estimates.
\begin{enumerate}
\item[(a)] 
$\displaystyle \Im \hat\tau_\gamma \geq \frac{1}{2},$
\item[(b)] 
$\displaystyle \log\Im\hat\tau_\gamma \leq \log \frac{d^2}{Nyk^2},\quad$ and
\item[(c)] 
$\displaystyle\log\Im\tilde\tau_\gamma \leq \log\Im\hat\tau_\gamma + \log 4.$
\end{enumerate}
\end{lemma}

\begin{proof}
We compute 
\begin{align*}
    \Im \hat\tau_\gamma 
    & = \frac{Ny}{d^2k^2}\cdot 
    \frac{1}{\left(\frac{Ny}{d^2}\right)^2 + \left(\frac{b}{d} - \frac{h}{k}\right)^2 } \\
    & = \frac{d^2}{Nyk^2} \cdot
    \frac{1}{1 + \frac{\left(\frac{b}{d} - \frac{h}{k}\right)^2}{\left(\frac{Ny}{d^2}\right)^2}} \\
    & = \frac{x}{1+t}
\end{align*}
and so
\[
    \log \Im \hat\tau_\gamma = \log \underbrace{\frac{d^2}{Nyk^2}}_{x}
    - \log \left(1 + \underbrace{\frac{\left(\frac{b}{d} - \frac{h}{k}\right)^2}{\left(\frac{Ny}{d^2}\right)^2}}_t\right).
\]
It follows that
\[
\log\Im\hat\tau_\gamma \leq \log \frac{d^2}{Nyk^2}.
\]

We also have $|\frac{b}{d} - \frac{h}{k}| \leq \frac{\sqrt{Ny}}{dk} $, so 

\[
0 \leq t \leq \frac{Ny}{d^2k^2} \cdot \frac{d^4}{N^2y^2} = \frac{d^2}{Nyk^2} = x.
\]

Furthermore, as $x \geq \frac{M^2}{k^2} \geq 1$, we also find
\[
\Im \hat\tau_\gamma \geq \frac{1}{2}.
\]

Finally, combining this with 
$-\frac{1}{2} < \Re \hat{\tau}_\gamma < \frac{1}{2}$ it follows that 
\[
\hat{\tau}_\gamma \in \CF \cup S\CF \cup ST^{-1}\CF \cup ST\CF,
\]
where
\[
S = \left(\begin{matrix}0 & -1 \\ 1 & 0 \end{matrix}\right)
\quad\text{and}\quad
T = \left(\begin{matrix}1 & 1 \\ 0 & 1 \end{matrix}\right)
\]
are the standard generators of $\SL_2(\BZ)$.
In particular, we find

\[
\log\Im\tilde\tau_\gamma \leq \log\Im\hat\tau_\gamma + \log 4.
\]
\end{proof}

Now we estimate the sum in (\ref{eq:target}) for those $\gamma\in C_N$ with $d \geq \sqrt{Ny}$. We note that $$2\log M \leq 2\log d - \log(Ny).$$

\begin{lemma}\label{large d}
Suppose $\tau = iy$ with $y\geq 1$ and $N\geq 1$. Then
$$\sum_{\substack{\gamma\in C_N\\ d_\gamma \geq \sqrt{Ny}}} \log\Im\tilde{\tau}_\gamma\leq \left( 4.75+ 3.5\log 2 + \frac{0.5 + \log 2}{2\sqrt{N}}\right)\psi(N).$$

\end{lemma}
\begin{proof}
    
Let us start with

\begin{align}
    \sum_{\substack{\gamma\in C_N\\ d_\gamma \geq \sqrt{Ny}}} \log\Im\tilde{\tau}_\gamma & \leq 
    \sum_{\substack{d|N \\ d \geq \sqrt{Ny}}} \; \sum_{\substack{0 \leq b < d\\ (b,r)=1}} \big(\log \Im \hat\tau_\gamma + \log4 \big)  \label{EqLargeD1}\\
    & \leq \left( \sum_{\substack{d|N \\ d \geq \sqrt{Ny}}} \; \sum_{\substack{0 \leq b < d\\ (b,r)=1}} \log \Im \hat\tau_\gamma \right) + \log(4)\psi(N), \nonumber
\end{align}
as the number of terms in the sum is bounded by $\#C_N = \psi(N).$ By Lemma \ref{tauhat},
\begin{align}
    \sum_{\substack{d|N \\ d \geq \sqrt{Ny}}} \; \sum_{\substack{0 \leq b < d\\ (b,r)=1}} \log \Im \hat\tau_\gamma
    & \leq \sum_d \sum_b \log \frac{d^2}{Nyk^2} \label{EqLargeD2} \\
    & = \sum_d \; \sum_{k=1}^M \; \underbrace{ \sum_{h=1}^k  \sum_{\frac{b}{d} \in I_M(\frac{h}{k})}}_{(i)} 
    \left[2\log\frac{d}{k} - \log(Ny)\right]. \nonumber
    \end{align}

Let us bound the number of terms in the sum $(i)$ above. By Lemma \ref{Farey}, the length of $I_M(\frac{h}{k})$ is bounded by $\frac{2}{(M+1)k}$, hence the number of terms in the inner sum is bounded by the number of integers $b$ with $(b,r)=1$ in $dI_M(\frac{h}{k})$, for fixed $k$ and $h=1,\ldots k$. For an interval of length $r$, we have $\varphi(r)$ integers coprime with $r$. Therefore, for fixed $d$ and $k$,
\begin{equation}\label{bound with+1}
\#\left\{ \frac{b}{d} \in I_M(\frac{h}{k}),\ 0 \leq b < d, \ (b,r)=1 \right\} \leq \varphi(r)\left\lceil\frac{1}{r}\frac{2d}{(M+1)k}\right\rceil \leq \varphi(r)\left(\frac{2d}{(M+1)kr} +1\right).
\end{equation}

Summing over $1\leq h \leq k$, we bound the number of terms in the sum $(i)$ by
\[
 \frac{2d\varphi(r)}{(M+1)r} + k\varphi(r).
\]
 
Hence we split the sum in Equation \eqref{EqLargeD2} in two:

\begin{gather}
\sum_{\substack{d|N\\ d \geq \sqrt{Ny}}} \; \sum_{\substack{0 \leq b < d \\ (b,r)=1}} \log \Im \hat\tau_\gamma \label{EqCeilingBound}\\
     \leq  \sum_d \sum_{k=1}^M \frac{2d\varphi(r)}{(M+1)r} \left[2\log\frac{d}{k} - \log(Ny)\right] + \sum_d \sum_{k=1}^M k\varphi(r) \left[2\log\frac{d}{k} - \log(Ny)\right] \nonumber
\end{gather}

We deal with the first sum in the right hand side of inequality \eqref{EqCeilingBound}.
    \begin{align}
&\sum_d \sum_{k=1}^M \frac{2d\varphi(r)}{(M+1)r} \left[2\log\frac{d}{k} - \log(Ny)\right] \label{EqLargeDShort}\\
    & = \sum_d \frac{2d\varphi(r)}{(M+1)r} \left[ 2M\log d - 2\sum_{k=1}^M\log k - M\log(Ny) \right] \nonumber \\
     &\leq \sum_d \frac{2d\varphi(r)}{r} \frac{M}{M+1} \left[ 2 + \log\frac{d^2}{M^2 Ny}-\frac{1}{6M(M+1)} - 2\frac{\log(\sqrt{2\pi M})}{M} \right] (ii) \nonumber \\
     &\leq \sum_d \frac{2d\varphi(r)}{r} \frac{M}{M+1} \left[ 2 + 2\log\frac{M+1}{M} -\frac{1}{6M(M+1)} - \frac{\log(2\pi M)}{M}  \right] \quad (iii) \nonumber \\
     &\leq \sum_{d|N} \frac{2d\varphi(r)}{r}2 \quad (iv)  \nonumber\\
     &= 4\psi(N), \nonumber
    \end{align}
\noindent where to reach ($ii$) we used 

$\displaystyle{\sum_{k=1}^{M} \log k= \log(M!)}$ and by \cite{Rob} we have for any integer $M\geq1$:
\[\sqrt{2\pi M}\left(\frac{M}{e}\right)^M e^{\frac{1}{12(M+1)}}\leq M!\leq  \sqrt{2\pi M}\left(\frac{M}{e}\right)^M e^{\frac{1}{12M}}\] so in particular
\[M \log M -M +\log\sqrt{2\pi M}+\frac{1}{12 (M+1)}\leq \sum_{k=1}^{M} \log k,\]
\noindent and in ($iii$) we used the fact that $M \leq \frac{d}{\sqrt{Ny}} \leq M +1 $ implies 
\[
1 \leq \frac{d^2}{M^2Ny} \leq \left(\frac{M+1}{M}\right)^2,
\]
and the inequality ($iv$) holds because for any $M\geq1$ we have
$$\frac{M}{M+1} \left[ 2 + 2\log\frac{M+1}{M} -\frac{1}{6M(M+1)} - \frac{\log(2\pi M)}{M}  \right]\leq2,$$
which can be verified through direct computation.

Let us bound the second sum in the right hand side of inequality \eqref{EqCeilingBound}.

\begin{align*}
&\sum_{\substack{d|N, \\ d \geq \sqrt{Ny}}} \sum_{k=1}^M k\varphi(r)\left(2 \log \frac{d}{\sqrt{Ny}} - 2\log k \right) \\
&= \sum_{\substack{d|N, \\ d \geq \sqrt{Ny}}} \left( \varphi(r)M(M+1) \log \frac{d}{\sqrt{Ny}} - 2\varphi(r)\sum_{k=1}^M k\log k \right).
\end{align*}

We bound $\displaystyle{-\sum_{k=1}^M k\log k}$ as follows.  By Abel's summation formula,
\[
\sum_{k=1}^M k\log k = \frac{M(M+1)}{2}\log M - \int_{1}^M \frac{\lfloor u \rfloor (\lfloor u \rfloor + 1)}{2}\frac{1}{u} du,
\]
hence, as $\frac{\lfloor u \rfloor}{u}\leq 1,$
\begin{gather*}
-2\sum_{k=1}^M k\log k = -M(M+1)\log M + \int_{1}^M \frac{\lfloor u \rfloor (\lfloor u \rfloor + 1)}{u} du \\\leq 
-M(M+1)\log M + \frac{M(M+1)}{2} - 1.
\end{gather*}
We set $\tilde{M} :=\frac{d}{\sqrt{Ny}}$, so $M \leq \tilde{M} \leq M+1$. Therefore,
\begin{align}
& \sum_{d|N, \; d \geq \sqrt{Ny}}\left( \varphi(r)M(M+1) \log \frac{d}{\sqrt{Ny}} - 2\varphi(r)\sum_{k=1}^M k\log k \right) \nonumber \\ 
& \leq \sum_{d|N, \; d \geq \sqrt{Ny}} \varphi(r) \left( M(M+1)\log \left( \frac{M+1}{M} \right)  + \frac{M(M+1)}{2} - 1 \right) \nonumber \\ 
& \leq \sum_{d|N, \; d \geq \sqrt{Ny}} \varphi(r) \left( \tilde{M}(\tilde{M}+1)\log 2  + \frac{\tilde{M}(\tilde{M}+1)}{2} \right) \nonumber\\
& = \sum_{d|N, \; d \geq \sqrt{Ny}} \varphi(r)\left( \left(\log 2 +\frac{1}{2}\right)\tilde{M}(\tilde{M} + 1)  \right)  \label{EqM(M+1)} \\
& \leq \left( \log 2 +\frac{1}{2} \right) \left(\sum_{d|N}\varphi(r)\frac{d^2}{Ny} + \sum_{d|N, \; d \geq \sqrt{N}}\varphi(r)\frac{d}{\sqrt{Ny}} \right) \label{EqFirstSumSecondSum} \\
&\leq \left( \log 2 +\frac{1}{2} \right)\left(\frac{3}{2} + \frac{1}{2\sqrt{N}}\right)\psi(N). \label{EqLastBound}
\end{align}
For the first sum in \eqref{EqFirstSumSecondSum}, remark that $a = \frac{N}{d}$ also runs through the divisors of $N$ and that $r = (d, a)$, hence $a \geq r$, and
\[
\sum_{d|N} \varphi(r)\frac{d^2}{Ny} = \frac{1}{y}\sum_{a|N} \varphi(r)\frac{N}{a^2} \leq \frac{1}{y}\sum_{a|N} \frac{\varphi(r)}{r}\frac{N}{a} = \frac{1}{y}\sum_{d|N} \frac{\varphi(r)}{r}d =  \frac{1}{y}\psi(N) \leq \psi(N).
\]
For the second sum in \eqref{EqFirstSumSecondSum},
\begin{equation}\label{inequalities}
\sum_{d|N, \; d\geq \sqrt{N}} \varphi(r) \frac{d}{\sqrt{Ny}} 
= \frac{1}{\sqrt{y}} \sum_{a|N, \; a \leq \sqrt{N}} \varphi(r)\frac{\sqrt{N}}{a} 
\leq \sqrt{N} \sum_{a\leq \sqrt{N}} \frac{\varphi(r)}{r} \leq \frac{\psi(N) + \varphi(\sqrt{N})}{2}.
\end{equation}
\noindent The last inequality in (\ref{inequalities}) comes from
\begin{align*}
\psi(N) &= \sum_{d > \sqrt{N}} \frac{\varphi(r)}{r}d + \sum_{d < \sqrt{N}} \frac{\varphi(r)}{r}d + \frac{\varphi(\sqrt{N})}{\sqrt{N}}\sqrt{N} \\
&= \sum_{d \leq \sqrt{N}} \frac{\varphi(r)}{r}\left( d + \frac{N}{d} \right) - \varphi(\sqrt{N}) \geq 2\sqrt{N} \sum_{d \leq \sqrt{N}} \frac{\varphi(r)}{r} -\varphi(\sqrt{N}) 
\end{align*}
as $(d + \frac{N}{d})\geq 2\sqrt{N}$ for any $1 \leq d \leq N$. We also set $\varphi(x) = 0$ if $x \not\in \mathbb{N}$. Notice also that $\frac{\varphi(\sqrt{N})}{\psi(N)} \leq \frac{1}{\sqrt{N}}$, because $\varphi(\sqrt{N}) \leq \sqrt{N}$ and $\psi(N) \geq N.$ Equation \eqref{EqLastBound} now follows.

This finishes the proof. 
\end{proof}

\begin{rem}
    We remark that we can obtain the slightly worse bound $(8 + 2\log 2)\psi(N)$ in Lemma \ref{large d} with a simpler argument. With the notations in \eqref{bound with+1}, it can be checked that the following inequality is true, $$\frac{2d}{(M+1)kr} \geq 1$$ for any $1 \leq k \leq M$. This implies that the second sum in \eqref{EqCeilingBound} is bounded by the first, so we could use the bound of Equation \eqref{EqLargeDShort} for both of them. 
\end{rem}
\subsection{Small $d$}

Now we consider the sum over $\gamma\in C_N$ with $d_\gamma < \sqrt{Ny}$. 

\begin{lemma}\label{small d}
Let $\tau = iy$ with $y\geq 1$ and $N\geq 1$. Then we have
$$\sum_{\substack{\gamma\in C_N \\ d_\gamma < \sqrt{Ny}}} \log\Im\tilde{\tau}_\gamma\leq \psi(N)\left(\frac{1}{e}+\log\Ima\tau\right).$$
\end{lemma}

\begin{proof}
In this case $\Im \tau_\gamma > 1$, so 
$\Im \tilde\tau_\gamma = \Im \tau_\gamma$. We write $a = \frac{N}{d}$ and compute

\begin{align*}
    \sum_{\substack{\gamma\in C_N \\ d_\gamma < \sqrt{Ny}}} \log\Im \tau_\gamma     
    & = \sum_{\substack{d|N \\ d < \sqrt{Ny}}} \; \sum_{\substack{b < d\\ (b,r)=1}} \log\frac{a y}{d} \\
    & = \sum_{\substack{d|N \\ d < \sqrt{Ny}}} \frac{d\varphi(r)}{r} \left[ \log\frac{a}{d} + \log y \right] \\
    & \leq \psi(N)\log y  + \sum_{\substack{d|N \\ d < \sqrt{Ny}}} \frac{d\varphi(r)}{r} \log\frac{a}{d}. 
\end{align*}

 The crude estimate $\log\frac{a}{d} \leq \frac{1}{e}\frac{a}{d}$ (which holds as $\frac{\log x}{x}$ has a maximum at $x=e$ for $x >0$) gives us
\begin{equation}
    \sum_{\substack{d|N \\ d < \sqrt{Ny}}} \frac{d\varphi(r)}{r} \log\frac{a}{d} 
    \leq \frac{1}{e}\sum_{a|N} \frac{a\varphi(r)}{r}  = \frac{1}{e}\psi(N).
\end{equation}
As $y=\Ima\tau$, this concludes the proof.
\end{proof}

\subsection{Final steps of the proof}

As shown in \cite{BP}, computations by Andrew Sutherland confirm Theorem \ref{thm:main1} for $N\leq 400$, so we may assume $N\geq 401$. In this case, the coefficient in Lemma \ref{large d} is
\[
4.75 + 3.5\log 2 + \frac{0.5+\log 2}{2\sqrt{N}} < 7.2059.
\]
Adding
the bounds in Lemma \ref{large d} and Lemma \ref{small d}, we obtain
\[
\sum_{\gamma\in C_N} \log\Im\tilde{\tau}_\gamma \leq \psi(N)(7.5737 + \log\Im\tau).
\]

We choose $\tau = iy$ such that $j(\tau) \in [1728, 3456]$, so $1\leq y < 1.2536$, for which we compute (using SageMath \cite{Sage}) 
\[
- \log\big[|\Delta(\tau)|(\Im\tau)^6\big] \leq 6.5296,
\]
and so in equation (\ref{eq:SN}) we have
\begin{align*}
S_N(\tau) & \leq  6\psi(N)\big[\log N - 2\lambda_N +  0.1878\big]
    + 6\sum_{\gamma\in C_N} \log\Im\tilde{\tau}_\gamma 
    - \psi(N)\log\big[|\Delta(\tau)|(\Im\tau)^6\big] \\
& \leq 6\psi(N)\big[\log N - 2\lambda_N +  9.0756\big].
\end{align*}

We add a classical interpolation lemma.

\begin{lemma}\label{interpolation}
Let $N\geq 1.$
For any real $L>1$, $$
h(\Phi_N) \leq \max_{L \leq j(\tau) \leq 2L} S_N(\tau) + \psi(N)\left(\frac{1+\log L}{L} + 4\log 2\right).
$$
\end{lemma}
\begin{proof}
    This is obtained in \cite{BP} equation (19), and comes from Lemma 10 in \cite{BrSu}.
\end{proof}

We can finally use Lemma \ref{interpolation} with $L=1728$ (corresponding to the smallest permissible value of $y$, which gives the best constants), and obtain
\begin{align*}
    h(\Phi_N) & \leq \max_{1728\leq j(\tau) \leq 3456} S_N(\tau) + \psi(N)\left(\frac{1+\log 1728}{1728} + 4\log 2\right) \\
    & \leq 6\psi(N) \big[\log N -2\lambda_N + 9.5387\big].
\end{align*}

This concludes the proof of the upper bound in Theorem \ref{thm:main1}.

\section{Proof of the lower bound in Theorem \ref{thm:main1}}\label{proof of lower}

We now turn to the lower bound in Theorem \ref{thm:main1}.
For any $\tau$ in the complex upper half plane, recall that the logarithmic Mahler measure of $\Phi_N(X,j(\tau))$ is equal to
\[
m\big(\Phi_N(X, j(\tau))\big) = S_N(\tau)=\sum_{\gamma\in{C_N}}\log\max\{1,\vert j(\tau_\gamma)\vert\}.
\]

We start with an upper bound that will be used later in the proof.

\begin{lemma}\label{mahler ll height}
    For every $N\geq 1$ we have
    \begin{enumerate}
        \item[(a)] $\displaystyle S_N(\tau) \leq 2\log(\psi(N)+1) + \psi(N)\log\max\{1, |j(\tau)|\} + h(\Phi_N).$ 
        \item[(b)] $\displaystyle S_N(\rho) \leq \log(\psi(N)+1) + h(\Phi_N)$, where $\rho = e^{\frac{i\pi}{3}}.$
    \end{enumerate}
\end{lemma}

\begin{proof}
     Write $\Phi_N(X,Y) = \sum_{k=0}^{\psi(N)}P_k(Y)X^k$, where each $P_k(Y)\in\BZ[Y]$ has degree $\leq \psi(N)$. Denoting $H(P_k) = e^{h(P_k)}$ the maximum absolute value of the coefficients of $P_k$, we have
    \[
    \big|P_k\big(j(\tau)\big)\big| \leq (\psi(N)+1)\max\{1, |j(\tau)|\}^{\psi(N)} H(P_k) \leq 
    (\psi(N)+1)\max\{1, |j(\tau)|\}^{\psi(N)} H(\Phi_N).
    \]
    Comparing the Mahler measure to the {\em length} of a polynomial \cite[Lemma 1.7]{BrZu}, we get
    \begin{align}\label{mahler upper bound}
    S_N(\tau) = m\big(\Phi_N(X,j(\tau))\big) & \leq \log\left[ \sum_{k=0}^{\psi(N)} \big|P_k\big(j(\tau)\big)\big| \right] \\ \nonumber
    & \leq \log \left[ (\psi(N)+1)^2 \max\{1, |j(\tau)|\}^{\psi(N)} H(\Phi_N) \right]. 
    \end{align}
    This proves part (a).

    Since $P_k(0)$ is the constant coefficient of $P_k(Y)$, we see that 
    \[
    \log |P_k(0)| \leq h(P_k) \leq h(\Phi_N).
    \]
    As $j(\rho)=0$, in this case (\ref{mahler upper bound}) gives
    \[
    S_N(\rho) \leq \log\left[\sum_{k=0}^{\psi(N)}|P_k(0)|\right] \leq \log(\psi(N)+1) + h(\Phi_N).
    \]
    Part (b) follows.
\end{proof}





To obtain a lower bound on $h(\Phi_N)$, it is thus enough to bound $S_N(\rho)$ from below, which is the goal of the next lemma.

\begin{lemma}\label{SNrho}
Let $\rho = e^{\frac{i\pi}{3}}$. Then for for any $N\geq1$, 
    $$S_N(\rho)\geq  6\psi(N)\left(\log N - 2\lambda_N -\frac{1}{6}\log\left\vert \frac{3^3}{(2\pi)^{24}}\Gamma\left(\frac{1}{3}\right)^{36}\right\vert-\frac{5.5335}{6}\right).$$
\end{lemma}

\begin{proof}
We bound the two sums in equation (\ref{eq:SN0}) for $S_N(\tau)$ from below, starting with Lemma \ref{logmaxbound} which gives us
 \begin{equation}
    \displaystyle{\sum_{\gamma\in C_N}} \log\max\{|\Delta(\tilde{\tau}_\gamma)|,|j(\tau_\gamma)\Delta(\tilde{\tau}_\gamma)|\} \geq -5.5335\psi(N).
\end{equation}
Also, for any $\gamma\in C_N$, we have $\Ima\tilde{\tau}_\gamma\geq \frac{\sqrt{3}}{2}$. We thus obtain




\begin{equation}\label{SNstep}
  S_N(\tau)\geq -5.5335\psi(N)  + 6\psi(N)(\log N - 2\lambda_N) + 6\psi(N)\log\frac{\sqrt{3}}{2} -\psi(N)\log\vert \Delta(\tau)(\Ima\tau)^6\vert.
\end{equation}

We will now specialize $\tau=\rho$. We obtain via Lemma \ref{deltarho} and equation (\ref{SNstep}):

\begin{equation}
  S_N(\rho)\geq -5.5335\psi(N) + 6\psi(N)(\log N - 2\lambda_N) -\psi(N)\log \left\vert \frac{3^3}{(2\pi)^{24}}\Gamma\left(\frac{1}{3}\right)^{36}\right\vert,
\end{equation}
which leads to the claim.
\end{proof}

By combining Lemma \ref{mahler ll height}(b) and Lemma \ref{SNrho}, we finally obtain 
$$h(\Phi_N)\geq -5.5335\psi(N) + 6\psi(N)(\log N - 2\lambda_N) -\psi(N)\log\left\vert \frac{3^3}{(2\pi)^{24}}\Gamma\left(\frac{1}{3}\right)^{36}\right\vert-\log(\psi(N)+1),$$ 
and 
$$\frac{1}{6}\left(\log\left\vert \frac{3^3}{(2\pi)^{24}}\Gamma\left(\frac{1}{3}\right)^{36}\right\vert 
+ \frac{\log(\psi(N)+1)}{\psi(N)}
 +5.5335\right)\leq 0.0351$$
 when $N \geq 401$. This proves the lower bound from Theorem \ref{thm:main1} in the case $N\geq 401$, whereas  the numerical computations by Andrew Sutherland (see \cite{BP}) show that the Theorem also holds when $N\leq 400$.

\section{Explicit Hecke points estimates}\label{Hecke}

So far, we have only obtained bounds on $\sum_{\gamma\in C_N}\log\Im\tilde\tau_\gamma$ for the special values $\tau = iy$. One can deduce a general bound from Theorem \ref{thm:main1}, which we record in the following result.

\begin{prop}\label{newprop}
    Let $\tau\in\CF$. Then
    \begin{enumerate}
        \item[(a)] $\displaystyle \max\big\{\log\frac{\sqrt{3}}{2}, \log\Im\tau - \log N + 2\lambda_N\big\}
        \leq \frac{1}{\psi(N)} \sum_{\gamma\in C_N}\log\Im\tilde\tau_\gamma \leq 10.832 + \log\Im\tau.$ 
        \item[(b)] If $\Im\tau \geq N$, then
        $\displaystyle\frac{1}{\psi(N)} \sum_{\gamma\in C_N}\log\Im\tilde\tau_\gamma = \log\Im\tau -\log N + 2\lambda_N.$
    \end{enumerate}
\end{prop}

\begin{proof}
    For each $\gamma\in C_N$ we have $\Im\tilde\tau_\gamma \geq \Im\tau_\gamma$, and also  $\Im\tilde\tau_\gamma \geq \frac{\sqrt{3}}{2}$. Thus
    \[
    \log\Im\tilde\tau_\gamma \geq \max\left\{\log\Im\tau_\gamma, \log\frac{\sqrt{3}}{2}\right\}.
    \]
    Now (\ref{eq:Aut2.3}) implies the lower bound in part~(a). 

    Furthermore, if $\Im\tau \geq N$, then $\Im\tilde\tau_\gamma = \Im\tau_\gamma$ for all $\gamma\in C_N$, and so (\ref{eq:Aut2.3}) implies part (b).

    We now prove the upper bound in part (a). 
    From Lemma \ref{mahler ll height}(a) we obtain
    \[
    S_N(\tau) \leq 2\log(\psi(N)+1) + \psi(N)\log\max\{1, |j(\tau)|\} + h(\Phi_N).
    \]  
    
    Next, we replace the left hand side by (\ref{eq:SN0}) and extract
    \begin{align}
        & \frac{1}{\psi(N)}\sum_{\gamma\in C_N}\log\Im\tilde\tau_\gamma \\ \nonumber  
        & \quad \leq  \left[\frac{1}{6\psi(N)}h(\Phi_N) - \log N + 2\lambda_N\right] \\ \nonumber
        & \quad +
        \frac{1}{6}\left[ 
        \log\max  \big\{|\Delta(\tau)|, |j(\tau)\Delta(\tau)|\big\} 
        -\frac{1}{\psi(N)}\sum_{\gamma\in C_N}\log\max  \big\{|\Delta(\tilde\tau_\gamma)|, |j(\tilde\tau_\gamma)\Delta(\tilde\tau_\gamma)|\big\}  
        \right] \\ \nonumber
        & \quad  + \log\Im\tau + \frac{\log(\psi(N) + 1)}{3\psi(N)}.
    \end{align}
    Now  Theorem \ref{thm:main1} and Lemma \ref{logmaxbound} give us
    \[
    \frac{1}{\psi(N)}\sum_{\gamma\in C_N}\log\Im\tilde\tau_\gamma \leq 
    9.5387 + \frac{1}{6}[1.1266 + 5.5335] + \frac{\log 3}{6} + \log\Im\tau.
    \]
    The result follows.
\end{proof}

We now prove Theorem \ref{thm:explicit Silverman}, which is an explicit version of Silverman's Theorem 5.1 page 417 of \cite{Sil}. 

\begin{proof}(of Theorem \ref{thm:explicit Silverman})
    Fix $N$ and $E$, and let $K$ be a sufficiently large number field that $E$ and every $E/C$ as well as the isogenies linking them are defined over $K$. 

    It follows from \cite[Prop. 2]{Sil} that only the infinite places contribute to the difference, so
    \begin{align}\label{height difference}
    & h_\infty(j_E)-\frac{1}{\psi(N)}\sum_{\substack{C\; cyclic \\ \# C=N}} h_\infty(j_{E/C}) \\ \nonumber
    = & \frac{1}{[K : \mathbb{Q}]}\sum_{\sigma : K \hookrightarrow \BC} \left[ \log\max\big\{1, |\sigma(j_E)|\big\} - 
    \frac{1}{\psi(N)}\sum_{\substack{C\; cyclic \\ \# C=N}} \log\max\big\{1, |\sigma(j_{E/C})|\big\} \right].
    \end{align}
    Notice that the Hecke sum in \cite{Sil} is over all subgroups $C \subset E$ of order $N$, not just the cyclic ones, but the argument in \cite[Prop. 2]{Sil} gives the same result in our situation.

    Let $\tau_\sigma\in\CF$ be such that $\sigma(j_E) = j(\tau_\sigma)$, then
    \[
    \sum_{\substack{C\; cyclic \\ \# C=N}} \log\max\big\{1, |\sigma(j_{E/C})|\big\} = S_N(\tau_\sigma)
    = m\big(\Phi_N(X,\sigma(j_E)\big)
    \]
    is the Mahler measure of $\Phi_N(X, j(\tau_\sigma))$. Now Lemma \ref{mahler ll height}(a) gives
    \[
    m\big(\Phi_N(X,\sigma(j_E))\big) \leq 
    2\log(\psi(N)+1) + \psi(N)\log\max\{1, |j_E|\} + h(\Phi_N).
    \]  
    Part (a) now follows from Theorem \ref{thm:main1} and the estimate 
    $\displaystyle\frac{\log(\psi(N)+1)}{\psi(N)} \leq \frac{\log 3}{2}$.

    To show part (b), we write $\tau = \tau_\sigma\in\CF$ and combine (\ref{eq:SN0}), Lemma \ref{logmaxbound} and Proposition \ref{newprop}:
    \begin{align} \label{mahler lower bound}
       & \log\max\big\{1, |\sigma(j_E)|\big\} - 
    \frac{1}{\psi(N)}\sum_{\substack{C\; cyclic \\ \# C=N}} \log\max\big\{1, |\sigma(j_{E/C})|\big\} \\ \nonumber
    & \quad = \log\max\big\{1, |j(\tau)|\big\} - \frac{1}{\psi(N)}S_N(\tau) \\ \nonumber
    & \quad = \left[
        \log\max\{|\Delta(\tau)|, |j(\tau)\Delta(\tau)|\} 
        - \frac{1}{\psi(N)}\sum_{\gamma\in C_N}\log\max\{|\Delta(\tilde\tau_\gamma)|, |j(\tilde\tau_\gamma)\Delta(\tilde\tau_\gamma)|\} 
        \right] \\ \nonumber
    & \qquad - \frac{6}{\psi(N)}\sum_{\gamma\in C_N}\log\Im\tilde\tau_\gamma 
    + 6\big[\log\Im\tau - \log N + 2\lambda_N\big] \\ \nonumber
    & \quad \leq [1.1266 + 5.5335]  \\ \nonumber
    & \qquad - 6\max\{\log\frac{\sqrt{3}}{2}, \log\Im\tau -\log N + 2\lambda_N\}  + 6\big[\log\Im\tau - \log N + 2\lambda_N\big] \\ \nonumber
    & \quad \leq 6.6601 + 6\min\{-\log\frac{\sqrt{3}}{2} + \log\Im\tau_\sigma - \log N + 2\lambda_N, 0 \}.
    \end{align}

    We now insert this into (\ref{height difference}) and invoke \cite[Lemma 2.6]{Paz}, which gives us
    \begin{align*}
    & h_\infty(j_E)-\frac{1}{\psi(N)}\sum_{\substack{C\; cyclic \\ \# C=N}} h_\infty(j_{E/C}) \\
    & \quad \leq 6.6601 + 6\min\big\{0, -\log\frac{\sqrt{3}}{2} + \log\big(1 + h_\infty(j_E)\big) + 1.94 -\log 2\pi - \log N + 2\lambda_N\big\}.
    \end{align*}
    This proves part (b) of Theorem \ref{thm:explicit Silverman}.

\end{proof}

\begin{rem}
    The inequalities (\ref{mahler lower bound}) and (\ref{mahler upper bound}) imply the following lower bound on the height of the specialised polynomial $\Phi_N(X,j)$, which can be seen as a measure of non-cancellation:
    \begin{align} \nonumber
        h(\Phi_N(X,j)) & \geq S_N(\tau) - \log(\psi(N)+1) \\ \nonumber
        & \geq \psi(N)\left[\log\max\{1, |j(\tau)|\} - 6.6601\right] - \log(\psi(N)+1) \\ 
        & \geq \psi(N)\left[\log\max\{1, |j(\tau)|\} - 7.2095  \right].
    \end{align}
\end{rem}

\begin{rem}
    If $N \leq \Im\tau_\sigma$ for every $\sigma : K \hookrightarrow \BC$, then the above proof, together with Proposition \ref{newprop}(b) gives
    \[
    \left|h_\infty(j_E)-\frac{1}{\psi(N)}\sum_{\substack{C\; cyclic \\ \# C=N}} h_\infty(j_{E/C})\right|
    \leq 6.6601.
    \]
\end{rem}

\begin{rem}\label{autissier}
Theorem \ref{thm:explicit Silverman} may be regarded as a ``Hecke-averaged'' version of \cite[Thm 1.1]{Paz}, with improved bounds. 

If we replace the Weil height of the $j$-invariant with the stable Faltings height (see Definition \ref{faltings}) of elliptic curves in Theorem \ref{thm:explicit Silverman}, Autissier \cite[Cor. 3.3]{Aut} obtained the even neater result:
\[
\frac{1}{\psi(N)} \sum_{\substack{C\; cyclic \\ \# C=N}} h_{\mathrm{Falt}}(E/C) 
= h_{\mathrm{Falt}}(E) + \frac{1}{2}\log N - \lambda_N.
\]
\end{rem}

\section{What is the size of $X_0(N)$?}\label{what?}

In this final section we give a proof of Theorem \ref{what is the size}. We start with the first item and recall the definition of the Faltings height of an abelian variety and of a curve.

\subsection{Faltings height and modular polynomials}
Let $A$ be a semi-stable abelian variety defined over a number field $k$, of dimension $g\geq 1$. Let $\pi\colon {\mathcal A}\longrightarrow \Spec(\mathcal{O}_k) $
be the N\'eron model of $A$ over $\Spec(\mathcal{O}_k)$, where $\mathcal{O}_k$ is the ring of integers of $k$. Let $\varepsilon\colon \Spec(\mathcal{O}_k)\longrightarrow {\mathcal A}$ be the zero section of $\pi$ and let
$\omega_{{\mathcal A}/\mathcal{O}_k}$ be the maximal exterior power of the sheaf of relative differentials
$$\omega_{{\mathcal A}/\mathcal{O}_k}:=\varepsilon^{\star}\Omega^g_{{\mathcal
A}/\mathcal{O}_k}.$$

For any archimedean place $v$ of $k$, let $\sigma$ be an embedding of $k$ in $\mathbb{C}$ associated to $v$. The associated line bundle
$$\omega_{{\mathcal A}/\mathcal{O}_k,\sigma}=\omega_{{\mathcal A}/\mathcal{O}_k}\otimes_{{\mathcal O}_k,\sigma}\mathbb{C}\simeq H^0({\mathcal
A}_{\sigma}(\mathbb{C}),\Omega^g_{{\mathcal A}_\sigma}(\mathbb{C}))\;$$
is equipped with a natural $L^2$-metric $\Vert.\Vert_{v}$ given by
$$\Vert s\Vert_{v}^2=\frac{i^{g^2}}{(2\pi)^{g}}\int_{{\mathcal
A}_{\sigma}(\mathbb{C})}s\wedge\overline{s}\;.$$

The ${\mathcal O}_k$-module $\omega_{{\mathcal A}/\mathcal{O}_k}$ is of rank $1$ and together with the hermitian norms
$\Vert.\Vert_{v}$ at infinity it defines an hermitian line bundle 
$\overline{\omega}_{{\mathcal A}/\mathcal{O}_k}=({\omega}_{{\mathcal A}/\mathcal{O}_k}, (\Vert .\Vert_v)_{v\in{M_k^\infty}})$ over $\mathcal{O}_k$. 

Recall that for any hermitian line
bundle $\overline{\mathcal L}$ over $\Spec(\mathcal{O}_k)$ the Arakelov degree of $\overline{\mathcal L}$ is defined as
$$\widehat{\degr}(\overline{\mathcal L})=\log\#\left({\mathcal L}/{s{\mathcal
O}}_k\right)-\sum_{v\in{M_{k}^{\infty}}}d_v\log\Vert
s\Vert_{v}\;,$$
where $s$ is any non zero section of $\mathcal L$. The resulting real number does not depend on the choice
of $s$ in view of the product formula on the number field $k$. 

The natural idea is then to consider
$\widehat{\degr}(\overline{\omega}_{{\mathcal A}/\mathcal{O}_k})$. This Arakelov degree of the metrized bundle $\overline{\omega}_{{\mathcal
A}/\mathcal{O}_k}$ will give a translate (by a term of the form $g c_0$ with $c_0$ an absolute constant) of the classical Faltings height. 

\begin{definition}\label{faltings}
  The stable height of $A$ is defined as
$$h_{\mathrm{Falt}}(A):=\frac{1}{[k:\mathbb{Q}]}\widehat{\degr}(\overline{\omega}_{{\mathcal
A}/\mathcal{O}_k})\;.$$
\end{definition}

In the same spirit, we can also define the Faltings height of a stable curve.

\begin{definition}
Let $k$ be a number field and $C/k$ a smooth algebraic curve defined over $k$, with semi-stable reduction and genus $g\geq1$. Let $p:C\rightarrow S$ be a semi-stable integral model of $C$ on $S=\Spec(\mathcal{O}_{k})$. The Faltings height of $C/k$ is the quantity
$$h_{\mathrm{Falt}}(C)=\frac{1}{[k:\mathbb{Q}]}\degrar (\det p_{*}\omega_{C/S}),$$
where the hermitian metrics are chosen as $\Vert\alpha\Vert_v^2=\frac{i^{g^2}}{(2\pi)^{g}}\int\alpha\wedge\overline{\alpha}$.
\end{definition}

This height is often referred to as the \textit{stable} height, as it is stable by extension of the base field $k$. The following proposition is well known to experts.

\begin{prop}\label{JacC}
Let $k$ be a number field and $C/k$ a smooth algebraic curve defined over $k$, with semi-stable reduction and genus $g\geq1$. Let $J_C$ denote the jacobian of $C$. Then we have
$$h_{\mathrm{Falt}}(J_C)=h_{\mathrm{Falt}}(C). $$
\end{prop}

\begin{proof}
See for instance Proposition 6.5 in \cite{Paz3}. 
\end{proof}

By specializing to $X_0(N)$, we get $h_{\mathrm{Falt}}(X_0(N))=h_{\mathrm{Falt}}(J_0(N))$. We now recall a result of Jorgenson and Kramer on the asymptotic of the Faltings height of the modular jacobian.
\begin{theorem} (Theorem 6.2 page 36 of \cite{JK})\label{JK}
    Let $N$ be square-free and coprime to 6. Let $g(N)$ be the dimension of the abelian variety $J_0(N)$. When $N$ tends to infinity, one has $$h_{\mathrm{Falt}}(J_0(N))=\frac{g(N)}{3}\log N +o(g(N)\log N).$$
\end{theorem}

We now need an estimate on the size of $g(N)$ as a function of $N$. This is done in the next lemma.

\begin{lemma}\label{genus}
    Let $N$ be square-free and coprime to 6. When $N$ tends to infinity, we have for any $\eps>0$ $$g(N)=\frac{N}{12}\prod_{p\vert N}\left(1+\frac{1}{p}\right)+O(\sigma(N)) = \frac{\psi(N)}{12}+O_{\eps}(N^{\eps}), $$
where $\displaystyle{\sigma(N) = \sum_{d \vert N}1.}$
For a general $N$,
$$g(N)= \frac{\psi(N)}{12} +O(\sqrt{N}\log\log(2N)).$$
\end{lemma}
\begin{proof}
    The dimension of $J_0(N)$ equals the genus of $X_0(N)$, which is given in Proposition 1.43 page 25 of \cite{Shi} by the formula, valid for $N$ coprime to 6, 
\begin{equation}
    g(N)=1+\frac{N}{12}\prod_{p\vert N}\left(1+\frac{1}{p}\right)-\frac{1}{4}\prod_{p\vert N}\left(1+\left(\frac{-1}{p}\right)\right)-\frac{1}{3}\prod_{p\vert N}\left(1+\left(\frac{-3}{p}\right)\right)-\frac{1}{2}\sum_{d\vert N}\varphi((d,\frac{N}{d})), \label{EqGenus}
\end{equation}
    where $\varphi$ is Euler's function and $\left(\frac{\cdot}{p}\right)$ is the quadratic residue symbol. In the general case, the formula has the same structure, with the products vanishing according to some divisibility conditions.
    
    Let us solve first the square-free case. One can check that the second and third products in the above expression either vanish or coincide with $\sigma(N)$ up to the corresponding constant factor in front of the product. With respect to the sum, if $N$ square-free then $(d, \frac{N}{d}) = 1$ for any $d \vert N$, and the sum equals $\sigma(N)$. The statement follows from the known growth rate of $\sigma(N) = O_{\eps}(N^{\eps})$ (see Theorem 315 from \cite{HW60}).
    
    In the general case, the products are still bounded by $\sigma(N)$. Define now
\[
\psitil(N) := \sum_{d | N} \varphi\left(\left(d, \frac{N}{d}\right)\right),
\]
where $(a,b)$ denotes the greatest common divisor of the integers $a$ and $b$, and $\varphi$ is Euler's totient function. Let us study this arithmetic function. Note that $\psitil$ is a multiplicative arithmetic function, i.e.\ $\psitil(ab) = \psitil(a)\psitil(b) \text{ if } (a,b)=1.$

For $p$ a prime number, $k \geq 1$ odd,
\begin{align*}
\psitil(p^k) &= \sum_{i=0}^{k} \varphi((p^i, p^{k-i})) = 2\sum_{i=0}^{\lfloor \frac{k}{2} \rfloor} \varphi(p^i) = 2( 1 + \sum_{i=1}^{\lfloor \frac{k}{2} \rfloor} (p^i - p^{i-1})) = 2p^{\lfloor \frac{k}{2} \rfloor} \\
&= 2p^{\frac{k}{2} - \frac{1}{2}} = \frac{2}{\sqrt{p}}\sqrt{p^{k}} \leq \left( 1 + \frac{1}{p} \right)\sqrt{p^{k}},
\end{align*}
since $\frac{2}{\sqrt{p}} < \left( 1 + \frac{1}{p} \right)$ for any prime.

Likewise, if $k \geq 1$ is even,
\begin{align*}
\psitil(p^k) &= \sum_{i=0}^{k} \varphi((p^i, p^{k-i})) = 2\sum_{i=0}^{\frac{k}{2} - 1} \varphi(p^i) + \varphi(p^{\frac{k}{2}}) = p^{\frac{k}{2} - 1} + p^{\frac{k}{2}} = \\
&= \left( 1 + \frac{1}{p} \right)\sqrt{p^{k}},
\end{align*}

Therefore, as $\psitil$ is multiplicative,
\[
\psitil(N) \leq \sqrt{N} \prod_{p | N}\left( 1 + \frac{1}{p} \right).
\]
As $\psi(N) = N\prod_{p | N} (1 + \frac{1}{p}),$
\begin{equation}
\psitil(N) \leq \frac{\psi(N)}{\sqrt{N}},
\label{EqBoundPsiPsitil}
\end{equation}
and it is known that $\psi(N) = O(N\log(\log(2N)))$ (see \cite{BD+96}[Lemme 2 (i)]). This finishes the proof.
\end{proof}
\begin{rem} It follows further from the proof of Lemma \ref{genus} that for any $\eps >0$, $\psitil(N) = O_{\eps}(N^{\frac{1}{2} + \eps})$ with explicit constant
\begin{equation}
C_{\eps} = \prod_{1 > p^{\eps} - \frac{1}{p}} p^{-\eps}\left( 1 + \frac{1}{p} \right). \label{EqCeps}
\end{equation}
We can therefore give an explicit (but worse) error term in the genus formula \eqref{EqGenus}. In particular, from \eqref{EqBoundPsiPsitil}, $\psi(N) \geq N$ and $\sigma(N) \leq 2\sqrt{N}$ we can deduce:
\begin{align*}
    \left| g(N) - \left( 1 + \frac{\psi(N)}{12}\right) \right| &\leq \frac{1}{2}C_{\eps}N^{\frac{1}{2} + \eps} + \frac{7}{12}\sigma(N) \leq \sqrt{N} \left(\frac{C_{\eps}}{2} N^{\eps} + \frac{7}{6} \right), \quad\text{and}\\
    \frac{ \left| g(N) - \left( 1 + \frac{\psi(N)}{12}\right) \right|}{\psi(N)} &\leq  \frac{1}{2}\frac{\psitil(N)}{\psi(N)} +\frac{7}{12}\frac{\sigma(N)}{\psi(N)} \leq \frac{5}{3}\frac{1}{\sqrt{N}}.
\end{align*}
It can also be shown, by inspecting how many primes verify the condition under the product in \eqref{EqCeps}, that the constant $C_{\eps}$ verifies:
\begin{itemize}
\item $C_{\eps} < 1$, for $\eps > 0.585$ (as the product is empty),
\item $C_{\eps} < 1.2527$ for $\eps > 0.26$ (as the product only has the prime $2$),
\item $C_{\eps} < 1.5788$ for $\eps > 0.132$ (as the product only has the primes $2$ and $3$).
\end{itemize}
\end{rem}

We can now conclude on the first item of Theorem \ref{what is the size}: by Proposition \ref{JacC}, $h_{\mathrm{Falt}}(X_0(N))=h_{\mathrm{Falt}}(J_0(N))$. By Theorem \ref{JK}, $h_{\mathrm{Falt}}(J_0(N))\sim \frac{g(N)}{3}\log N$ when $N$ tends to infinity and is square-free, coprime to 6. By Lemma \ref{genus}, $g(N)\sim \frac{\psi(N)}{12}$. Use the Corollary page 390 of \cite{Coh} which gives $h(\Phi_N)\sim 6\psi(N)\log N$ to conclude that 
\begin{equation}\label{first item}
h_{\mathrm{Falt}}(X_0(N))\sim \frac{1}{6^3} h(\Phi_N).
\end{equation}

\subsection{Hecke correspondences and modular polynomials}
Let us move to the second item of Theorem \ref{what is the size}. In \cite{Aut}, Autissier uses a morphism $i_N: X_0(N)\longrightarrow \mathbb{P}^1\times\mathbb{P}^1$, which for two elliptic curves $E_1, E_2$ and a cyclic isogeny $\alpha: E_1\to E_2$ is defined by $i_N((E_1,E_2, \alpha))=(j(E_1), j(E_2))$. He denotes by $T_N$ the image of $X_0(N)$ by $i_N$, and by $\hat{\mathcal{L}}$ a natural metrized lined bundle on $\mathbb{P}^1\times\mathbb{P}^1$. Theorem 3.2 page 427 of \cite{Aut} gives $$h_{\hat{\mathcal{L}}}(T_N)=12\psi(N)(\log N-2\lambda_N+4\kappa_1),$$ where $\kappa_1=12\zeta'(-1)-\log\pi-\frac{1}{2}$, which implies that for any $N\geq1$, $\vert h_{\hat{\mathcal{L}}}(T_N)- 2 h(\Phi_N)\vert$ is bounded by a quantity linear in $\psi(N)$, which in turn implies, as the main term is of order of magnitude bigger than $\psi(N)$, the fact that when $N$ tends to infinity $$h_{\hat{\mathcal{L}}}(T_N)\sim 2 h(\Phi_N).$$

\subsection{Heegner points and modular polynomials}
The third item in Theorem \ref{what is the size} comes from an asymptotic estimate computed in \cite{Paz2} and heavily based on the Gross-Zagier computations \cite{GZ}. Corollaire 1 page 164 in \cite{Paz2} provides us, when $N$ tends to infinity and satisfies the Heegner conditions (there are infinitely many such $N$ for each fixed discriminant $D_k$), with $$\hat{h}_{J_0(N)}(c_{D_k})\sim \frac{3h_k u_k}{g(N)}h_{\mathrm{Falt}}(J_0(N)),$$ where $\hat{h}_{J_0(N)}$ is the N\'eron-Tate height on the jacobian $J_0(N)$ as defined in \cite{GZ}. Use $h_{\mathrm{Falt}}(J_0(N))=h_{\mathrm{Falt}}(X_0(N))$ and (\ref{first item}) to obtain this third item. 
This concludes the proof of Theorem \ref{what is the size}.

\subsection{Arakelov canonical sheaf of $X_0(N)$}
The fourth item in Theorem \ref{what is the size} comes from the following asymptotic estimate, first computed in Th\'eor\`eme 1.1 page 646 of \cite{MU} in the case where $N$ is coprime to 6 and square-free, and recently generalised to any N coprime to $6$ in Theorem 1.1 of \cite{DM}:
\begin{equation}
    \overline{\omega}^2\sim 3g(N)\log N.
\end{equation}
As we have $g(N)\sim \frac{\psi(N)}{12}$ by Lemma \ref{genus} and by Corollary page 390 of \cite{Coh} we have $h(\Phi_N)\sim 6\psi(N)\log N$, hence we get the result.


\end{document}